\documentclass{article}
\usepackage[T2A]{fontenc}
\usepackage[utf8]{inputenc}
\usepackage[english,russian]{babel}
\usepackage[tbtags]{amsmath}
\usepackage{amsfonts,amssymb,mathrsfs,amscd}
\usepackage[hyper]{msb-a}

\makeatletter
\gdef\No{{\select@language{russian}\textnumero}}
\makeatother

\JournalName{}
\numberwithin{equation}{section}

\def \Z {{\mathbf {Z}}}

\def \T {{\mathbf  T}}
\def \R {{\mathbf  R}}
\def\uu{\bigsqcup}
\def\eps{\varepsilon}

\begin{document}

\title{Полиномиальная жесткость и  спектр  сидоновских автоморфизмов}
\author[V.\,V.~Ryzhikov]{В.\,В.~Рыжиков}
\address{Московский государственный университет}
\email{vryzh@mail.ru}

\date{02.11.2023}
\udk{517.987}

\maketitle

\begin{fulltext}

\begin{abstract}{Предъявлен континуум  спектрально дизъюнктных сидоновских автоморфизмов, тензорный  квадрат которых изоморфен планарному сдвигу. Спектры таких автоморфизмов не обладают групповым свойством.  Для выявления сингулярности спектров использована  полиномиальная  жесткость операторов,  связанная с понятием линейного детерминизма  по  Колмогорову.   В классе перемешивающих гауссовских и пуассоновских надстроек над  
сидоновскими автоморфизмами реализованы   новые наборы спектральных кратностей.

Библиография: 12 названий.}
\end{abstract}
\markright{Полиномиальная жесткость и  спектр }

\section{Введение} 
\bf Обобщение свойства  жесткости оператора. \rm  Пусть $T:H\to H$ -- ограниченный оператор в нормированном пространстве $H$, обладающий следующим свойством:
для некоторой  последовательности $n_j\to +\infty$ и  плотного в   $H$ подпространства $F$    найдется последовательность полиномов $ P_j(T)=\sum_{n, |n|>n_j} a_n T^n$  такая, что для всякого  $f\in F$
$$ \|P_j(T)f -f\|\to 0.$$
Говорим, что такой  оператор $T$ обладает  $P_j$-жесткостью.
В этом случае также говорим, что оператор обладает обобщенной (полиномиальной)  жесткостью.  
Напомним, что оператор $T$ называется  жестким, если для некоторой последовательности $n_i\to+\infty$ выполнено $T^{n_j}\to_w I$. Для унитарного оператора $T$ это свойство влечет за собой  сингулярность его спектральной меры. 

 Полиномиальная жесткость операторов в  гильбертовых пространствах тесно связано с линейным детерминизмом  по Колмогорову, т.е. возможностью восстанавливать вектор в виде линейной  комбинации его далеких сдвигов (см. \cite{K}, а также \cite{H},\cite{O} о развитии этой тематики).
 
 \bf Полиномиальное перемешивание. \rm    Последовательность  полиномов вида  $$ P_j(T)=\sum_{n, |n|\leq n_j} a_{j,n} T^n, \ a_{j,n}\geq 0,\ n_j\to+\infty,$$ 
будем называть \it перемешивающей,  \rm если для некоторого плотного  подпространства $V\subset H$ для любого $v\in V$ выполнено $\|P_j(T)v\|\to 0$. Терминология взята из    эргодической теории:  свойство перемешивания автоморфизма $T$ вероятностного пространства эквивалентно тому, что степени индуцировнного им унитарного оператора слабо сходятся к 0 в пространстве, ортогональном константам.  Мы  рассматриваем случай, когда  $H=L_2(\mu)$, где мера $\mu$ сигма-конечна, поэтому ненулевые константы не входят в $L_2(\mu)$.

\bf Дизъюнктность операторов. \rm Спектральные меры  жесткого оператора $T$ вдоль последовательности $n_j$ и  перемешивающего оператора $\tilde T$ вдоль последовательности $n_j$ взаимно сингулярны.  Это свойство эквивалентно  их дизъюнктности: для ограниченного
 оператора  $J$ из $TJ=J\tilde T$  вытекает $J=0$.  Действительно,  переходя к слабому пределу в $T^{n_j}J=J\tilde T^{n_j}$, получим $J=0$.  Модификация этого соображения позволяет различать  спектры  перемешивающих операторов  $T$ и $\tilde T$.  Если найдется  последовательность  полиномов  $P_j$ такая, что $P_j(T)$ является жесткой  последовательностью операторов, а $P_j(\tilde T)$  перемешивающая, то   $T$ и $\tilde T$ дизъюнктны.

 \bf  Отсутствие группового свойства спектра. \rm В 60-е годы А.Н. Колмогоров высказал гипотезу    о групповом свойстве спектра автоморфизмов, которая, получив отрицательное решение  (см.  \cite{KS}),  стимулировала появление новых подходов в  изучении спектров  эргодических действий. 
В предлагаемой статье мы  рассмотрим   так называемые простые сидоновские  автоморфизмы $T$,  сохраняющие сигма-конечную меру.  Будет показано, что  их спектры сингулярны, а сверточные квадраты спектров при дополнительных условиях  эквивалентны лебеговской мере. Таким образом,  мы предъявим новые примеры автоморфизмов $T$, для которых свертка $\sigma_T\ast\sigma_T$  взаимно сингулярна с их  спектральной мерой $\sigma_T$. 
 \bf  Тензорные корни планарного сдвига. \rm  В связи с задачей  Ж.-П. Тувено о неизоморфных тензорных корнях тензорного квадрата автоморфизма (см. \cite{18}) будет  установлено, что  ненулевой сдвиг $S:\R^2\to\R^2$,  рассматриваемый как преобразование, сохраняющее  плоскую меру  Лебега,  имеет спектрально неизоморфные  тензорные корни.  Мы покажем, что  для   континуума спектрально дизъюнктных сидоновских автоморфизмов $T:\R\to\R$ их  прозведения  $T\times T$ сопряжены со сдвигом  $S$ некоторым автоморфизмом $\Phi_T$ (т.е.  $T\times T=\Phi_T S\Phi^{-1}_T$). Сидоновские конструкции рассматривались  в работе   \cite{20} в  качестве примеров   метрически неизоморфных тензорных корней преобразования $S$. 

\bf Приложения для пуассоновских и гауссовских динамических систем.  \rm 
Автоморфизму $S$  стандартного пространства с сигма-конечной мерой  отвечают две классические конструкции автоморфизмов вероятностного пространства: гауссовский автоморфизм $G(S)$ и пуассоновская надстройка $P(S)$, хорошо известные не только в эргодической теории,
 но и в теории представлений групп (см., например, \cite{N}). 
 В работе \cite{23} были  реализованы всевозможные спектральные кратности  вида  $M\cup\{\infty\}$ для неперемешивающих гауссовских систем и пуассоновских надстроек. Результаты настоящей статьи приводят к аналогичной реализации в классе перемешивающих систем.

 \section{ Вспомогательные утверждения. Конструкции автоморфизмов ранга один} Сингулярность спектра  сидоновских автоморфизмов (определение 
которых дано в \S 3)  можно получить, следуя   методу и результатам  работы \cite{DE}. Отметим, что для других классов автоморфизмов сингулярность  спектра доказывалась в \cite{B}, \cite{KR}. Мы не применяем   обобщенные  произведения Рисса, а  вместо произведений используем  суммы полиномов, которые обеспечивают  обобщенную  жесткость. Для простого сидоновского автоморфизма  будет  доказана  дизъюнктность его степеней и как следствие     сингулярность  его спектра.  В основе предлагаемого подхода лежат следующие     утверждения.

\vspace{2mm}
\bf  Лемма 2.1. \it Пусть для унитарных  операторов $S$ и $T$  последовательность  $P_j(T)$ является жесткой, 
а последовательность  $P_j(S^\ast)$   перемешивающая. Тогда
операторы $S$ и $T$ дизъюнктны: для всякого ограниченного оператора $J$, удовлетворяющего 
условию  $SJ=JT$,  вытекает $J=0$. \rm

\vspace{2mm}
Лемма   является частным случаем следующего  утверждения.

\vspace{2mm}
\bf  Теорема 2.2. \it Пусть    для унитарных операторов $S$ и $T$  на пространстве $H$ найдутся плотные в $H$  подпространства $F,G$ такие, что последовательность   $P_j(T)f$ сходится  к $P(T)f$ при $f\in F$, а $P_j(S^\ast)g$ сходится к $P'(S^\ast)g$ при $g\in G$. Если полиномы  $P$ и $P'$ различны, а  оператор $T$  имеет непрерывный спектр,
то  $T$ и $S$  дизъюнктны. \rm

\vspace{2mm}
Доказательство.   Пусть $SJ=JT$. 
Имеем $$(P(T)f,J^\ast g)=\lim_j (JP_j(T)f,g)= \lim_j(Jf, P_j(S^\ast)g)=$$
$$=(Jf,P'(S^\ast) g)= (P'(S)Jf, g)= (P'(T)f,J^\ast g).$$
Следовательно, $$( P(T)f-P'(T)f,J^\ast g)=0,  \ \ f\in F,\, g\in G.$$
Так как $T$ имеет непрерывный спектр, векторы $P(T)f-P'(T)f$, $f\in F$, плотны в $H$,
получим $J^\ast=0$.

\vspace{2mm}
\bf  Лемма 2.3. \it Пусть для унитарного оператора $T$ 
последовательность $P_j(T)$ жесткая, а  для всех $p>1$ последовательности  $P_j(T^p)$ перемешивающие, причем  все полиномы имеют вещественные  коэффициенты.
Тогда  спектр оператора $T$ сингулярен. \rm

\vspace{2mm}
Доказательство. Для  оператора $T$  выполнено $\|P_j(T^p)v\|=\|P_j(T^{-p})v\|$ (что очевидно, когда у многочлена $P_j$ вещественные коэффициенты). Поэтому из леммы 2.1. вытекает, что  $T$ дизъюнктен со степенями $T^p$, $p>1$. При наличии  у оператора $T$ абсолютно непрерывной компоненты $\sigma_{1}$ последняя  для всех больших  $p$  не  может быть взаимно сингулярна с   абсолютно непрерывной компонентой $\sigma_p$ оператора $T^p$.
Действительно,  пусть $\sigma_1$ задается плотностью Радона-Никодима
$\rho_1(z)$,  соответствующая плотность $\rho_p(z)$ для оператора $T^p$  
является суммой $\sum_{w, w^p=z} \rho_1(w)/p$ (мера   $\sigma_p$  -- образ  $\sigma_1$  при отображении $z\to z^p$). Мера тех точек $z$, для которых $\rho_p(z)=0$, стремится к 0 при $p\to\infty$. Следовательно, для больших $p$ мера множества точек $z$ таких, что $\rho_1(z)>0$ и $\rho_p(z)>0$, положительна. Лемма доказана.

\vspace{2mm}
\bf  Конструкции    автоморфизмов ранга один. \rm
Напомним необходимое для дальнейшего определение автоморфизма ранга один.
Фиксируем натуральное  число $h_1\geq 1$ (высота башни на этапе $j=1$), последовательность  $r_j\to\infty$ ( параметр $r_j$ -- число колонн, на которые виртуально разрезается башня этапа $j$)   и последовательность целочисленных векторов (параметров надстроек)   
$$ \bar s_j=(s_j(1), s_j(2),\dots, s_j(r_j-1),s_j(r_j)).$$  
Ниже дано описание конструкции сохраняющего меру преобразования, которое полностью определено параметрами $h_1$, $r_j$ и $\bar s_j$.

На шаге $j=1$ задан  полуинтервал $E_1$. На  шаге $j$  определена  
 система   непересекающихся полуинтервалов 
$$E_j, TE_j, T^2E_j,\dots, T^{h_j-1}E_j,$$
причем на $E_j, TE_j, \dots, T^{h_j-2}E_j$
пребразование $T$ является параллельным переносом. Такой набор   полуинтервалов  называется башней этапа $j$, их объединение обозначается через $X_j$ и тоже называется башней.

Представим   $E_j$ как дизъюнктное объединение  $r_j$ полуинтервалов 
$$E_j^1,E_j^2,E_j^3,\dots E_j^{r_j}$$ одинаковой меры (длины).  
Для  каждого $i=1,2,\dots, r_j$ рассмотрим так называемую колонну  
$$E_j^i, TE_j^i ,T^2 E_j^i,\dots, T^{h_j-1}E_j^i.$$
К каждой  колонне с номером $i$  добавим  $s_j(i)$  непересекающихся полуинтервалов (этажей)  длины, равной длине интервала $E_j^i$.
Полученные  наборы интервалов  при  фиксированных $i$,$j$ называем надстроенными колоннами  $X_{i,j}$. Отметим, что при фиксированном $j$ по построению колонны $X_{i,j}$   не пересекаются. Используя параллельный пренос интервалов, преобразование $T$ теперь  доопределим так, чтобы колонны $X_{i,j}$ имели вид 
$$E_j^i, TE_j^i, T^2 E_j^i,\dots, T^{h_j-1}E_j^i, T^{h_j}E_j^i, T^{h_j+1}E_j^i, \dots, T^{h_j+s_j(i)-1}E_j^i,$$
  а   верхние этажи $T^{h_j+s_j(i)-1}E_j^i$  колонн  $X_{i,j}$ ($i<r_j$) преобразование $T$ параллельным переносом  отображало в нижние
этажи колонн $X_{i+1,j}$: 
$$T^{h_j+s_j(i)}E_j^i = E_j^{i+1}, \ 0<i<r_j.$$ 
Положив $E_{j+1}= E^1_j$, замечаем, что все указанные этажи надстроенных колонн в новых обозначениях имеют вид 
$$E_{j+1}, TE_{j+1}, T^2 E_{j+1},\dots, T^{h_{j+1}-1}E_{j+1},$$
 образуя башню  этапа $j+1$ высоты  
 $$ h_{j+1} =h_jr_j +\sum_{i=1}^{r_j}s_j(i).$$

Таким образом,  преобразование $T$ теперь определено на всех
этажах башни этапа $j+1$,  кроме последнего. Полученное частичное определение преобразования $T$ на этапе $j$ сохраняется на всех последующих этапах. В результате   получаем пространство  $X=\cup_j X_j$ и обратимое преобразование $T:X\to X$, сохраняющее  стандартную меру Лебега $\mu$  на $X$.  
 Если  
$$\sum_j {(s_j(1)+s_j(2)+\dots+s_j(r_j))}/{h_jr_j }=\infty,$$
то мера пространства $Х$  будет бесконечной, именно этот случай  рассматривается в статье.
Преобразования ранга один эргодичны, в случае бесконечной меры фазового пространства они имеют  непрерывный спектр.  Более того, хорошо известно, что линейные комбинации индикаторов  интервалов, фигурирующих  в описании конструкций  $T$, являются циклическими векторами для оператора $T$ 
(преобразование и индуцированный им унитарный оператор в статье обозначаются одинаково).  Разнообразные результаты об автоморфизмах  ранга один обсуждаются в \cite{20}.

\section{Сидоновские автоморфизмы и диссипативность  тензорного квадрата}  Пусть конструкция $T$ ранга один обладает следующим свойством: для всех  больших $j$ \it пересечение 
$X_j\cap T^mX_j$ при   $h_{j}<m\leq h_{j+1}$  может содержаться 
только в одной  из колонн  $X_{i,j}$ башни $X_j$. \rm Такая конструкция 
называется  \it сидоновской \rm в силу следующей аналогии: подмножество $M$  натурального ряда называется сидоновским, если $|M\cap M+n|\leq 1$ для всех $n>0$.
Заметим, что  для такого $T$ и   множества $A$, состоящего из этажей башни этапа $j_0$,
 выполняется $$\mu(A\cap T^mA)\leq \mu(A)/r_j$$
для всех   $m\in [h_j, h_{j+1}],$ $j\geq j_0$.  
Следовательно,   если  $r_j\to\infty$, сидоновская конструкция 
обладает перемешиванием.  Напомним, что автоморфизм  $T$ пространства с сигма-конечной мерой   \it перемешивающий\rm: если для всех множеств  $A,B$ конечной меры выполнено 
$$\mu(A\cap T^mB)\to 0.$$

Позже мы сосредоточим наше внимание на так  называемых простых сидоновских конструкциях, параметры которых обеспечивают указанное свойство пересечений $X_j\cap T^mX_j$. Для них в случае, когда $m$ сравнимо с $h_{j+1}$,
$X_j\cap T^mX_j\subset X_{1,j}$.  

\bf Пример сидоновской конструкции.  \rm Положим    $$s_j(i)=10^ih_j - h_j, \ 1\leq i\leq r_j.  $$
Заметим, что пересечение $X_j\cap T^mX_j$ непусто при $h_{j}<m\leq h_{j+1}$  только в случае,  когда число $m$  сравнимо 
с $(10^k-10^i)h_j$, более того, будет  выполнено 
$$T^m X_j\cap X_j= T^mX_{i,j}\cap  X_{k,j}, \ 1\leq i< k<r_j.\eqno (1)$$ 
Здесь наблюдается аналогия с таким сидоновским множеством как $M=\{10^q: 
1\leq q\leq r\}$, для которого $(M\cap M+m)\neq \phi$ только при  
$m=10^k-10^i\leq 1$.  Следует сказать, что из-за  надстроек на шаге $j+1$   мы имеем  $X_j\cap T^mX_j=\phi$
при  $0.12h_{j+1}<m< 9 h_{j+1}$.  Непустые  пересечения $T^m X_j\cap X_j$ при $h_{j}<m\leq h_{j+1}$ удовлетворяют условию   $(1)$, а это означает, что наша   конструкция является сидоновской.  Забегая вперед, отметим, что  конструкция, например,  с параметрами 
$$s_j(i)=(10+j)^ih_j, \ 1\leq i\leq r_j,$$ 
 обладает свойством: для всех больших $j$ и фиксированного $p$ из $T^m X_j\cap X_j\neq \phi$ при $m>h_j$ вытекает  $T^{pm} X_j\cap X_j= \phi$. Это свойство, как мы покажем позже,  влечет за собой  спектральную дизъюнктность степеней автоморфизма $T$.

Если для сидоновской конструкции последовательность $r_j$ растет достаточно быстро, то произведение $T\times T$ оказывается диссипативным (по этой причине оператор $T\otimes T$ имеет лебеговский спектр).

\vspace{2mm}
\bf Теорема 3.1.   \it Пусть для сидоновской конструкции $T$ выполнено
$$\sum_{i=1}^{\infty }{r_j^{-1}} \ <\infty,$$ 
Тогда $T\times T$ является диссипативным преобразованием с блуждающим множеством бесконечной меры. 
\rm

\vspace{2mm}
Доказательство.  
Наша цель -- установить при  $S=T\times T$ равенство 
$$X\times X=\uu_{i\in\Z} S^iW$$ 
для некоторого (блуждающего) множества $W$ бесконечной меры. Предварительно мы найдем блуждающее множество конечной меры.
Фиксируем натуральное $p$. Всюду далее используем обозначение  $\bar \mu=\mu\otimes\mu$. 
Мера множества возвращения в $C_p=X_p\times X_p$  за время  от  $h_j$ до $h_{j+1}$ при $j\geq p$  не превосходит $\bar \mu(C_p)/r_j$. Действительно, условие $S^iC_p\cap C_p\neq\phi$ равносильно  $T^iX_p\cap X_p\neq\phi$, причем 
$$S^iC_p\cap C_p\subset (T^hX_p\cap X_p)\times (T^hX_p\cap X_p).$$
При $h_j\leq h < h_{j+1}$ для сидоновской конструкции 
пересечение $T^hX_p\cap X_p$ лежит в одной из колонн $ X_{i,j}$, следовательно, 
$$S^h C_p\cap C_p \subset X_{i,j}\times X_{i,j}$$
для некоторого $i=i(h)$, $1<i\leq r_j$.
Так как $$ \bar \mu\left(\bigcup_{i=1}^{r_j}X_{i,j}\times X_{i,j}\right)=\frac {\bar \mu(C_p)} {r_j},$$  получаем, что в $C_p$ с момента 
времени $h_p$  возвращается множество меры, не больше  
$\sum_{j=p}^{\infty }{\bar\mu(C_p)}/ {r_j}.$
Так как ряд $\sum_{i=1}^{\infty }r_j^{-1}$ сходится, из сказанного следует, что мера  максимального блуждающего множества $W$ не меньше, чем  
$(1-\eps_p)\bar\mu(C_p)\to\infty$, $\eps_p\to 0$. Поясним это.
В силу сказанного выше мы имеем
$$\limsup_{N\to\infty} \bar \mu\left(S^{N+1}C_p \setminus\bigcup_{i=1}^N S^iC_p\right)=
\bar\mu(C_p)\to\infty,   \ p\to\infty.$$
Если   мера максимального блуждающего множества $W$ конечна, то 
для всякого $Y$, $ \bar \mu(Y)<\infty$,  верно, что 
$$\limsup_{N\to\infty} \bar \mu\left(S^{N+1}Y \setminus\bigcup_{i=1}^N S^iY\right)\leq 
\bar\mu(W).$$
Значит, в нашем случае мера $W$ бесконечна.  Теорема доказана.
 
\vspace{2mm}
В работе \cite{LS} даны примеры автоморфизмов $T$, обладающих консервативной (и даже эргодической) степенью $ (T\times\dots\times T)$ произвольного порядка такой, что произведение $ (T\times\dots\times T)\times T$ диссипативно.

\bf Замечания о консервативности и  спектре тензорных степеней.  \rm   Произведение $T\times T$ консервативно для сидоновских автоморфизмов, удовлетворяющих условиям 
$ r_j\to\infty$,    $ \sum_{i=1}^{\infty }r_j^{-1} = \infty.$  
В этом классе автоморфизмов  \it  для всякого натурального $m$ найдется  такой  $T$, что    $T^{\otimes m}$ имеет  сингулярный спектр, а для   некоторого $n>m$  степень $T^{\otimes n}$ диссипативна.  \rm 

С помощью таких  автоморфизмов можно получать новые спектральные типы  эргодических действий на вероятностном пространстве.
Например, \it для всякого $n>1$ найдется перемешивающая пуассоновская надстройка и перемешивающий гауссовский автоморфизм, для которых  спектральный 
тип имеет вид $\sigma +\sigma^{\ast 2}+\dots + \sigma^{\ast n}+\lambda,$
где сверточная степень $\sigma^{\ast n}$ сингулярна, 
а $\lambda$ -- мера Лебега.  \rm

Этим результатам  автор планирует посвятить отдельную работу.

 \section{Полиномиальная жесткость   и перемешивание для простых сидоновских автоморфизмов}
 Пусть $\psi>3$ и на  всех этапах $j$, начиная с некоторого,  для параметров конструкции ранга один выполнено 
$$ s_j(1)>\psi h_j,  \ \ s_j(i+1)> \psi s_j(i),   \ 1\leq i <r_j.$$
Несложно проверить, что такая конструкция  сидоновская. 
Будем называть    \it простой сидоновской \rm конструкцию,   параметры  
которой удовлетворяют условию
$$ h_j\ll  s_j(1)\ll  s_j(2)\ll \dots\ll  s_j(r_j-1)\ll s_j(r_j), \ \ r_j\to\infty,$$ 
где выражения вида $a_j\ll b_j$  обозначают, что $b_j>\psi(j)a_j$ 
для некоторой фиксированной последовательности $\psi(j)\to +\infty$.
Отметим свойство простой сидоновской конструкции, важное для дальнейшего.
 Из $$T^nX_1\cap X_1\neq\phi, \ \ h_j< n\leq 2s_j(r_j-1),$$ вытекает,  что для некоторого $i<r_j$ 
указанное пересечение $T^nX_1\cap X_1$ содержится в  колонне $X_{i+1,j}$  и выполнено 
$$T^nX_1\cap X_1= T^n(X_1\cap X_{i,j})\cap (X_1\cap X_{i+1,j}).$$
В этом случае $$ s_j(i)< n< s_j(i)+2h_j. $$

Так как  $h_j\ll s_j(i)\ll s_j(i+1)$,  множество
$T^{pn}(X_1\cap X_{i,j})$  окажется в надстройке над колонной 
$X_{i+1,j}$ и не будет пересекаться с $X_1\subset X_j$. Более того, пересечение множества $X_1$ со  всеми колоннами с номерами не больше $i$
под действием $T^{pn}$ окажется в этой же надстройке.
Сказанное  (при фиксированном $p$ для всех достаточно больших $j$)  вытекает из неравенства
$$ s_j(1)+h_j + s_j(2)+h_j +\dots +s_j(i)+h_j  < pn< s_j(i+1). $$
Таким образом, выполнено   $$T^{pn}X_1\cap X_1=\phi.$$
Как мы увидим, это свойство будет причиной того, что
для подходящих полиномов  $P_j$ мы получим  одновременно свойство
$P_j(T)$-жесткости  и $P_j(T^p)$-перемешивания. 
Для простых сидоновских конструкций $T$ мы найдем такую последовательность полиномов
$ P_j$, что для некоторого циклического вектора $f$  будет выполняться 
$ \|P_j(T)f -f\|\to 0$,
но при этом $ \|P_j(T^p)f\|\to 0$ для всякого $p>1$. Из  этого вытекает, что спектр $T$
сингулярный.

\vspace{2mm}
\bf  Теорема 4.1. \it Пусть $T$ -- простой  сидоновский автоморфизм с параметрами $s_j(i)$, $r_j$.   Положим $k(1,j)=0$  и   $$ k(i,j)=(i-1)h_j +s_j(1)+s_j(2)+\dots +s_j(i-1),\, 1<i\leq r_j,$$ 
$$Q_j(T)=\frac 1{r_j}\sum_{1\leq i\neq i'\leq r_j}  T^{k(i,j)-k(i',j)},$$
$$P_j=(Q_{j+1}+Q_{j+2}+\dots+Q_{j+m(j)})/m(j)$$
для некоторой последовательности  $m(j)\to\infty.$

Тогда последовательность $P_j(T)$  жесткая, а  последовательность 
$P_j(T^p)$ при $p>1$ является перемешивающей. \rm

\vspace{2mm}
\bf  Следствие. \it   Простой сидоновский автоморфизм $T$  имеет сингулярный спектр. \rm

\vspace{2mm}

Доказательство.  Обозначим через  $f$   индикатор  башни $X_1$. Докажем, что 

(i)\ \ $P_j(T^p)f\to 0$ при  $p>1$; 

     (ii) \ \    $P_j(T)f\to f$.
\\
Так как спектр $T$ непрерывен, для автоморфизма ранга  один такая функция  $f$ является циклическим вектором (хорошо известный факт).  Таким образом, установив  (i) и (ii) для   $f$, мы докажем (i) и (ii) для множества функций, образующих плотное подпространство в $L_2(X,\mu)$. Тем самым теорема будет доказана.

\bf Доказательство (i).  \rm   Пересечение  $T^{pk(i,j)-pk(i',j)}X_1\cap X_1$ при $1\leq i\neq i'\leq r_j$  для всех достаточно больших $j$ пусто. Это вытекает из того, что 
 в случае непустого пересечения 
$pk(i,j)-pk(i',j)/k(m,j)\approx 1$ для некоторого $m$.  Но это невозможно при $p>1$ силу условий
$$h_j\ll k(1,j)\ll k(2,j) \dots \ll k(r_j,j),$$
которые выполняются очевидным образом.
Аналогично убеждаемся в том, что 
$$T^{pk(i,j)-pk(i',j)-pk(m,j)+pk(m',j)}X_1\cap X_1\neq \phi$$ только при  $i=i'$ и $m=m'$, следовательно, пересечение  совпадает с $X_1$.

Ниже под знаком   суммы всегда подразумевается, что $1\leq i,i',m,m'\leq r_j$.  Например,  
запись $\sum_{ i\neq i'} $ означает   $\sum_{1\leq i\neq i'\leq r_j}$.
 
Для $f=\chi_{X_1}$ получаем 
$$\|Q_j(T^p)f\|^2=
\frac 1{r_j^2}\sum_{i\neq i'}  \sum_ {m\neq m'}
\left(T^{pk(i,j)-pk(i',j)-pk(m,j)+pk(m',j)}f\, , \, f\right)= $$

$$=\frac 1{r_j^2}\sum_{i=m\neq i'=m'} \mu(X_1)= (r_j-1)\|f\|^2/r_j.$$
Заметим, что множество $T^{pk(i,j)-pk(i',j)}X_1$ лежит в $X_{j+1}\setminus X_j$ (для всех больших $j$). Из этого вытекает, что носители функций $Q_{j+l}(T^p)f$, $l=1,2,\dots, m(j)$,
не пересекаются. Значит, эти фукции попарно ортогональны. Таким образом, 
$$\|P_j(T^p)f\|< \|f\|/\sqrt{m(j)}\to 0,  \ m(j)\to\infty.$$

\bf Доказательство (ii). \rm Специфика сидоновских конструкций выражается в том, что отображение $T^{k(i,j)-k(i',j)}$ отправляет $i'$-тую колонну $X_{i',j}$ в  $i$-тую колонну $X_{i,j}$, а образы  остальных колонн при  $i\neq i'$ оказываются в $X_{j+1}\setminus X_j$. 

Нам важно знать, что 
 $$\mu(T^{k(i,j)-k(i',j)-k(m,j)+k(m',j)}X_1\cap X_1)\in \{0,\, \mu(X_1)/r_j,\,\mu(X_1)\}$$
Случай ненулевой меры возможен, когда 
$k(i,j)-k(i',j)-k(m,j)+k(m',j)$ отличается от некоторой разности вида $k(n,j)-k(n',j),$ не более, чем  на $h_j$. Вспоминая, что   $h_j\ll k(1,j)\ll k(2,j) \dots \ll k(r_j,j),$
получим, что указанное условие выполняется только в случае, если $n,n'$ совпадают с некоторыми из чисел $i,i'm,m'$.  
Но это влечет за собой  совпадение двух чисел из набора $i,i'm,m'$. Пусть, например,  $i=m$,
а $i'\neq m'$.
Замечаем, что  пересечение
$T^{k(i,j)-k(i',j)}X_1 \cap T^{k(i,j)-k(m',j)}X_1$
   содержит  $X_1\cap X_{i,j}$, что 
проверяется непосредственно.
Тогда множество 
$$T^{k(i,j)-k(i',j)- k(i,j)+k(m',j)}X_1\cap X_1 =T^{k(m',j)-k(i',j)}X_1\cap X_1$$
содержит  $X_1\cap X_{m',j}$ и на самом деле совпадает
с ним, так как такое пересечение может содержаться лишь в одной из колонн $X_{i'',j}$.
Таким образом, мы убедились в том, что мера множества 
$T^{k(i,j)-k(i',j)}X_1 \cap T^{k(i,j)-k(m',j)}X_1$ при $i=m$,
 $i'\neq m'$  равна $\mu(X_1)/r_j$. Если $i=m$, $i'= m'$, очевидно, что мера пересечения  равна   $\mu(X_1)$.
 
Вычислим  $\|Q_j(T)f\|$ для функции   $f=\chi_{X_1}$. Имеем   
$$\|Q_j(T)f\|^2=
\frac 1{r_j^2}\sum_{i\neq i', m\neq m'}
\left(T^{k(i,j)-k(i',j)-k(m,j)+k(m',j)}f\, , \, f\right)= $$
$$=\frac 1{r_j^2}\left(\sum_{i=m\neq   i'= m'} \mu(X_1)\,+\,
 \sum_{i=m,   i'\neq m'} \frac {\mu(X_1)}{r_j}\, +\, 
 \sum_{i\neq m,    i'=m'}  \frac {\mu(X_1)}{r_j}\right)=$$
$$= 3(r_j-1)\|f\|^2/r_j.$$
Положим 
$$\Delta_j=Q_jf -\left(1-\frac 1{r_j}\right)f.$$

Для различных  $j$ функции $\Delta_j$ попарно ортогональны, так как  их носители не пересекаются  (они лежат в
$X_{j+1}\setminus X_j$). 

Таким образом, 
 $$\|\Delta_j\|^2=2\left(1-\frac 1{r_j}\right)\|f\|^2, $$ 
$$\left\|P_j(T)f-\left(1-\frac 1{r_j}\right)f\right\|^2  <\frac {2\|f\|^2} {m(j)}.$$
C учетом
 $r_j\to\infty$ получаем
$$\|P_j(T)f-f\| \ \to\ 0,  \ m(j)\to\infty.$$
Из доказанных сходимостей (i), (ii) вытекает утверждение теоремы, благодаря которой  в силу лемм 2.1 и 2.3  получаем сингулярность спектра автоморфизма $T$.

\vspace{2mm} \bf 
Замечание. \rm Если положить 
$$Q_j^+(T)=\frac 1{2r_j}\sum_{1\leq i'< i\leq r_j}  T^{k(i,j)-k(i',j)},$$  то для 
$P_j^+=(Q_{j+1}^++\dots+Q_{j+m(j)}^+)/m(j)$   снова получим приближение  
$$\|P_j^+(T)f-f\| \ \to\ 0,$$
 хотя для этого потребуется дополнительная аргументация.
Тем самым реализуется  полиномиальная жесткость при условии, что степени и коэффициенты полиномов $P_j^+$ положительны. В связи со сказанным  возникает  вопрос: \it влечет ли 
свойство $P_j(T)$-жесткости  при этом условии  сингулярность спектра автоморфизма $T$? \rm

\section{Дизъюнктные сидоновские автоморфизмы, тензорные корни, пуассоновские надстройки}
\bf Теорема 5.1.  \it Найдется континуальное семейство сидоновских автоморфизмов со взаимно сингулярными однократными спектрами, причем их всевозможные  тензорные произведения диссипативны (имеют лебеговский спектр).  
Найдется континуум  спектрально неизоморфных тензорных корней диссипативного преобразования, изоморфного планарному сдвигу. \rm

\vspace{2mm}
\bf Замечание. \rm  Если указано континуальное семейство автоморфизмов с взаимно сингулярными спектральными мерами, то  среди них  найдется не более, чем счетное  множество  автоморфизмов  с абсолютно непрерывной компонентой в спектре. Почти все (в смысле мощности) автоморфизмы семейства должны  иметь сингулярный спектр.

Доказательство теоремы. Фиксируем некоторую простую сидоновскую  конструкцию $T$. Напомним, что для ее  параметров  выполнено  $r_j\to\infty$, $ h_j\ll s_j(1)\ll  \dots \ll s_j(r_j)$), где $h_j$ -- высота  башни на этапе построения с номером $j$. Потребуем также выполнения условия 
$\sum_{i=1}^{\infty }{r_j^{-1}} \ <\infty,$ которое влечет за собой диссипативность 
$T\times T.$

\bf Фиктивные этапы с $\bf r_j=1$. \rm Приступим к   определению семейства   сидоновских автоморфизмов $\tilde T$,  для которых будет установлена попарная дизъюнктность их спектров.
  В определении конструкции ранга один требуется естественное условие  $r_j\geq 2$. Сейчас для удобства  мы  разрешим  случай $r_j=1$, предполагая, что  $r_j\geq 2$ имеет место  для бесконечного множества этапов $j$.   Этапы $j$ с $r_j=1$ будем называть  фиктивными. Мы просто увеличиваем параметр $s_{j-1}(r_j)$  этапа $j-1$  и называем это этапом $j$. Параметр увеличивается так, чтобы  для  высот башен конструкции $\tilde T$  выполнялось $\tilde h_{j}=h_{j}$. 
 
Обозначим $G_n = \{n^2,n^2+1,\dots, n^2+2n\}$. Пусть $\gamma$  -- некоторое бесконечное подмножество натурального ряда. Зададим параметры конструкции $\tilde T=T_\gamma$:   

если  $j\in G_n$, $n\notin \gamma$,  положим $\tilde r_j=1$  и  
$\tilde s_j(1)=h_{j+1}-h_j$; 

если $j\in G_n$, $n\in \gamma$, то все параметры конструкции $T_\gamma$ на этапе $j$ совпадают с параметрами конструкции $T$. 

 Рассмотрим полиномы, фигурирующие в теореме 4.1, 
положив  $k(i,j)=(i-1)h_j +s_j(1)+s_j(2)+\dots +s_j(i)$, $s_j(0)=0$ и  определяя 
$$Q_j(T)=\frac 1{r_j}\sum_{1\leq i\neq i'\leq r_j}  T^{k(i,j)-k(i',j)}$$
и 
$$P_j=(Q_{j^2+1}+Q_{j^2+2}+\dots+Q_{j^2+j})/ j.$$

Повторяя рассуждения в доказательстве теоремы 4.1,  получим, что  при $j_n\to\infty$, $j_n\in \gamma$ последовательность $P_{j_n}(T_\gamma)$   жесткая.  При $j_n\to\infty$, $j_n\notin \gamma$ последовательность $P_{j_n}(T_{\gamma'})$ перемешивающая.

Фиксируем континуальное семейство $\Gamma$  бесконечных подмножеств натурального ряда таких, что все пересечения различных предствителей семейства имеют конечную мощность.\ footnote{Такое семейство строится, например,  как аналог континуального семейства  "просек  в лесу    $\Z^2$", когда пересечение различных "просек" есть конечное множество "деревьев"  внутри параллелограмма.}  Тогда семейство автоморфизмов $T_\gamma$, $\gamma\in \Gamma$, будет искомым.  Если $\gamma\neq\gamma'$,  то при $j_n\in \gamma$ последовательность $P_{j_n}(T_\gamma)$ жесткая, но 
последовательность $P_{j_n}(T_{\gamma'})$ перемешивающая.
Тогда из  леммы 2.1 вытекает  дизъюнктность автоморфизмов 
 $T_\gamma$ и $ T_{\gamma'}$.

Остается заметить, что    произведения  $T_\gamma\times T_{\gamma'}$
 являются диссипативными.  Если $\gamma\neq  \gamma'$,
то для больших $k$ для башни $X_1$ выполнено:
$T_\gamma^kX_1\cap X_1\neq\phi$ влечет за собой 
$T_{\gamma'}^kX_1\cap X_1=\phi$ и тем самым диссипативность
произведения $T_\gamma\times T_{\gamma'}$.   Причины диссипативности  произведения 
$T_\gamma\times T_{\gamma}$  те же самые, что и для произведения $T\times T$  в доказательстве  теоремы 3.1.  Во всех случаях мера блуждающего множества равна бесконечности, поэтому 
все  произведения $T_\gamma\times T_{\gamma'}$ будут изоморфны планарному сдвигу. Теорема доказана.

\vspace{2mm}
\bf Спектральные кратности перемешивающих гауссовских и  пуассоновских надстроек.  \rm 
В  \cite{23} предъявлены спектральные  кратности вида  $M\cup\{\infty\}$ для неперемешивающих гауссовских  $G(S)$ и пуассоновских надстроек   $P(S)$ над автоморфизмом $S$ пространства с сигма-конечной мерой. 
Из теоремы 5.1 и  известного факта о том, что $G(S)$ и $P(S)$  унитарно эквивалентны оператору 
$$\exp(S)=1\oplus S\oplus S^{\odot 2}\oplus S^{\odot 3}\oplus \dots,$$ 
вытекает следующее утверждение.

\vspace{2mm}
\bf Теорема 5.2. \it Для всякого подмножества $M$ натурального ряда
найдется перемешивающая пуассоновская надстройка и перемешивающий гауссовский автоморфизм, у которых набор спектральных кратностей есть $M\cup\{\infty\}$. При этом бесконечную кратность
имеет лебеговская компонента. \rm

\vspace{2mm}
Действительно, $M\cup\{\infty\}$ является набором спектральных кратностей для оператора 
$$P=\exp\left(\bigoplus_{m\in M} \bigoplus_{i=1}^m  T_m\right), $$
при условии, что сингулярные спектры операторов $T_m$ однократны и попарно дизъюнктны, а все произведения вида $T_m\otimes T_{m'}$ имеют лебеговский спектр.  

Пуассоновская  надстройка   над  прямыми суммами $\bigoplus_{m\in M} \bigoplus_{i=1}^m  T_m$  является динамической реализацией оператора $P$ (аналогично для гауссовского  автоморфизма). Подходящий набор  автоморфизмов $T_m$ найдется благодаря  теореме  5.1.

\section{Замечания о детерминированности  сдвига в $l_p(\Z)$.}
Свойство $P_j(S)$-жесткости  является инвариантом оператора 
$S$ и по этой причине может представлять самостоятельный интерес.  В качестве примера покажем, как устанавливается это свойство  для  оператора  сдвига $S$ в банаховом пространстве  $l_p(\Z)$ при $p>2$.

\vspace{2mm}
\bf Утверждение 6.1.  \it При  $p>2$  сдвиг $S$  в $l_p(\Z)$  обладает обобщенной жесткостью.
\rm

\vspace{2mm}
Доказательство. Определим двустороннюю последовательность $v:\Z\to \R$  такую, чтобы величины  $v(i)$ принимали ненулевые значения    только при 
$i=-10, - 10^2, \dots, -10^n, \dots$. Наложим дополнительное условие:  все ненулевые значения  
$v(i)$ имеют вид $2^{-j}$, где $j=1,2,\dots$,  причем значение $2^{-j}$ среди  всех значений $v(i)$  встречается $q(j)= [2^{cj}]$ раз, где число $c$ фиксированно, причем  $2<c<p$.  
Для удобства будем считать, что ненулевые значения координат вектора $v$ убывают слева направо. Тогда найдется последовательность $h_j\to\infty$ такая, что  для  $i\in [h_j, h_{j+1})$ имеет место $v(i)=2^{-j}$ или $v(i)=0$, а для всех  $i\notin[h_j, h_{j+1})$  выполнено  $v(i)\neq 2^{-j}$. 
Так как  $c<p$, получим
$$\|v\|_p^p = \sum_{j=1}^\infty   q(j)2^{-pj}<\sum   2^{(c-p)j}<\infty.$$
Обозначим через $e_i$  орт, для которого $e_i(i)=1$, а остальные координаты нулевые.
Положим
$$Q_j(S)=\sum_{n, 10^n\in [h_j, h_{j+1})}   S^{10^n},$$
тогда
$$Q_j(S)v=\sum_{n, 10^n\in [h_j, h_{j+1})} 2^{-j} S^{10^n}v = q(j)2^{-j}e_0 + \Delta_j,$$
где    $\Delta_j$ --  сумма $q(j)$ векторов,  нормы которых меньше $\|v\|$, а их носители не пересекаются. Имеем
$$\|\Delta_j\|_p< q(j)^{\frac 1 p}  \|v\|_p \leq   2^{cj/p} \|v\|_p.$$
 С учетом того, что $q(j)= [2^{cj}]$  и  $2<c<p$,  получим 
$$ q(j)^{-1}2^{j}\|\Delta_j\|_p\ \leq   q(j)^{-1} 2^{j-cj/p}
 \|v\|_p\to 0.$$ Тогда 
$$R_j(S)v:=q(j)^{-1}2^{j}Q_j(S)v = e_0 +  q(j)^{-1}2^{j}\Delta_j\to e_0.$$

Так как  $S^iR_j(S)v$ сходится к $ e_i$,  всякий финитный вектор $\sum_i a_ie_i=\sum_i a_iS^ie_0$  можно приблизить векторами  $R_j(S)\sum_i a_iS^iv$. Финитные векторы плотны в $l_p$, следовательно, некоторая подходящая последовательность полиномов $P_j(S)$ является жесткой, т. е. $P_j(S)v\to v$, причем    коэффициенты  полиномов $P_{j}$ неотрицательны,
а  степени $S$, входящие в линейную комбинацию $P_{j}(S)$ с положительными коэффициентами, стремятся  к $+\infty$ при возрастании  $j$.  Утверждение доказано. 

\vspace{2mm}
\bf Сдвиг в $\bf l_2(\Z)$. \rm   Оператор сдвига в $l_2(\Z)$ изоморфен  оператору умножения $U$ на $z$ в пространстве $L_2(\T,\lambda)$, $|z|=1$, где $\T$ -- единичная окружность в комплексной плоскости с мерой Лебега $\lambda$. Установим полиномиальную жесткость для оператора $U$.
Из теоремы Колмогорова-Сегё (см., например, обзор \cite{O}) 
 вытекает существование ограниченной положительной на $\T$ весовой функции $w(z)$ такой, что для всякого $j$ система $\{z^n: n>j \}$  полна в $L_2(\T,w)$.  В силу полноты таких систем найдутся  полиномы $P_j$, линейные комбинации функций $z^n$  при $n>j$ такие, что  $P_j\to 1$  в $L_2(\T,w)$. Это означает, что 
$$\int_\T |P_j(z)-1|^2w(z)d\lambda(z) \to 0,$$
что влечет за собой 
$$\int_\T |P_j(z)w-w|^2d\lambda(z)= \int_\T |P_j(z)-1|^2w(z)^2d\lambda(z) \to 0.$$
Тем самым мы получили  $P_jw\to w$  в $L_2(\T,\lambda)$.   Так как  функция $w$ не обращается в 0, она является циклическим вектором для оператора  $U$. Возьмем в  качестве плотного пространства $F$  линейную оболочку векторов $z^nw$, $n\in \Z$, тогда для всякого $f\in F$ будет выпонено   $P_jf\to f$  в $L_2(\T,\lambda)$. Мы установили жесткость последовательности полиномов $P_j$.
В этом примере в отличие от предыдущих  не утверждается, что коэффициенты полиномов $P_j$ положительны.

\vspace{4mm}
\bf Благодарности. \rm Автор  выражает признательность  Ж.-П. Тувено и  участникам семинара  механико-математического  факультета МГУ по теории функций действительного переменного  за  полезные обсуждения. Автор  благодарит рецензента за замечания.

\end{fulltext}

\newpage
V.V. Ryzhikov

\vspace{3mm}
Polynomial rigidity and  spectrum of Sidon automorphisms.

\vspace{3mm}
A continuum of spectrally disjoint Sidon automorphisms is presented whose tensor squares are isomorphic to a planar shift. The spectra of such automorphisms do not have the group property. To identify the singularity of spectra, we use polynomial rigidity of operators associated with the concept of Kolmogorov linear determinism. In the class of mixing Gaussian and Poisson superstructures over
new sets of spectral multiplicities are realized by Sidon automorphisms.

\vspace{3mm}
Keywords:    Sidon automorphisms,  spectrum and  disjointness of transformations,  tensor roots and products,  polynomial rigidity and  mixing.


\begin{thebibliography}{99}
\bibitem{K}А.Н. Колмогоров, Интерполирование и экстраполирование стационарных случайных последовательностей, Изв. АН СССР. Сер. матем., 5:1 (1941), 1-40; A. N. Kolmogorov, Interpolation and extrapolation of stationary random sequences, Izv. Academy of Sciences of the USSR. Ser. Mat., 5:1 (1941), 1-40;


\bibitem{H} H. Helson.	M\' ethodes complexes et m\' ethodes de Hilbert en analyse de Fourier. Facult\' e des Sciences d'Orsay, 1967

\bibitem{O}А.М. Олевский, Представление функций экспонентами с положительными частотами, УМН, 59:1(355) (2004),  169-178; 
A.M. Olevskii, Representation of functions by exponentials with positive frequencies, Russian Math. Surveys, 59:1 (2004), 171-180

\bibitem{KS}	А.Б. Каток, А.М. Степин, Аппроксимации в эргодической теории, УМН, 22:5(137) (1967),  81-106; A.B. Katok, A.M. Stepin, Approximations in ergodic theory, Russian Math. Surveys, 22:5 (1967), 77-102

\bibitem{18} В.В. Рыжиков, Задача Тувено об изоморфизме тензорных степеней эргодических потоков, Матем. заметки, 104:6 (2018),  912-917;
V.V. Ryzhikov, Thouvenot's isomorphism problem for tensor powers of ergodic flows, Math. Notes, 104:6 (2018), 900-904 


\bibitem{20} 	В.В. Рыжиков, О сохраняющих меру преобразованиях ранга один, Тр. ММО, 81:2 (2020),  281-318; V.V. Ryzhikov, Measure-preserving rank one transformations, Trans. Moscow Math. Soc., 81:2 (2020), 229-259

\bibitem{N}	Ю. А. Неретин, Категории симметрий и бесконечномерные группы, УРСС, М., 1998



\bibitem{23}В.В. Рыжиков, Спектры самоподобных эргодических действий, Матем. заметки, 113:2 (2023),  273-282; V.V. Ryzhikov, Spectra of Self-Similar Ergodic Actions, Math. Notes, 113:2 (2023), 274-281


\bibitem{DE}A. H. Dooley, S. J. Eigen,  A family of generalized Riesz products, Can. J. Math. Vol. 48 (2), 1996 pp. 302-315 

 \bibitem{B}	J. Bourgain, On the spectral type of Ornstein's class one transformation, Israel J. Math., 84 (1993), 53-63


\bibitem{KR}
I. Klemes, K. Reinhold, Rank one transformations with singular spectral type, Israel J. Math., 98 (1997), 1-14 


\bibitem{LS}I. Loh,  C.E. Silva, Strict doubly ergodic infinite transformations, 
Dyn. Syst. 32, No. 4, 519-543 (2017). 


\end{thebibliography}
\end{document}